\documentclass[conference]{IEEEtran}
 
\usepackage{cite}
\ifCLASSINFOpdf
   \usepackage[pdftex]{graphicx}
\else
   \usepackage[dvips]{graphicx}
\fi
\usepackage{amsmath, amsfonts, amssymb}
\usepackage{url}

\usepackage[]{fancyhdr} % 
\newcommand{\changefont}{\fontsize{9}{9}\selectfont}
\usepackage{optidef, bm}
\usepackage{xcolor}

\IEEEoverridecommandlockouts
\begin{document}

\title{Uncertainty-Aware Methods for Leveraging Water Pumping Flexibility for Power Networks}

\author{\IEEEauthorblockN{Anna Stuhlmacher \& Johanna L. Mathieu}
\IEEEauthorblockA{Electrical Engineering and Computer Science\\
University of Michigan, Ann Arbor, MI, USA\\
\{akstuhl, jlmath\} @ umich.edu}
\thanks{This work was supported by NSF Grant 1845093.}}


%Telephone: (800) 555--1212, Fax: (888) 555--1212}

% <-this % stops a space


%\lhead{11TH BULK POWER SYSTEMS DYNAMICS AND CONTROL SYMPOSIUM, JULY 25-30, 2022, BANFF, CANADA}
%\rhead{1}

%\fontfamily{phv}\fontseries{b}\fontsize{9}{11}\selectfont

\maketitle
\thispagestyle{fancy}
\pagestyle{fancy}

%\thispagestyle{fancy}
%\pagestyle{fancy}

\begin{abstract}
%145/150 words
Recent work has demonstrated that water supply pumps in the drinking water distribution network can be leveraged to provide flexibility to the power network, but existing approaches are computationally demanding and/or overly conservative.  In this paper, we develop a computationally tractable probabilistic approach to schedule and control water pumping to provide voltage support to the power distribution network subject to power and water distribution network constraints under power demand uncertainty. Building upon robust and chance-constrained reformulation approaches,  we analytically reformulate the probabilistic problem into a deterministic one and solve for the scheduled pump operation and the control policy parameters that adjust the pumps based on the power demand forecast error realizations. In a case study, we compare our proposed approach to an adjustable robust method and investigate the performance in terms of computation time, cost, and empirical violation probabilities. We find that our proposed approach is computationally tractable and is less conservative than the robust approach, indicating that our formulation would be scalable to larger networks. 
\end{abstract}

\begin{IEEEkeywords}
distribution networks, flexible loads, uncertainty management, voltage support, water networks
%The author shall provide up to 5 keywords (in alphabetical order) to help identify the major topics of the paper. The thesaurus of IEEE indexing keywords is posted at \url{http://www.ieee.org/organizations/pubs/ani_prod/keywrd98.txt.}
\end{IEEEkeywords}


% \fi
\IEEEpeerreviewmaketitle

\section{Introduction} \label{section: intro}
There is growing interest in the coordinated operation and control of interdependent critical infrastructure systems, such as electric power grids, water networks, and natural gas networks. Leveraging the interconnection between multiple systems can lead to increased reliability, reduced operational costs, and improved sustainability over all systems \cite{DallAnese-Magazine, Mancarella2016}. However, there are also several challenges associated with the integrated operation of multiple systems. First, there is a need for more communications and measurements between traditionally independent system operators. Second, the problem complexity and dimension increase as we consider (nonconvex) network models from multiple systems and the propagation of uncertainty across interconnected systems. This can make the problems very difficult to solve. Approximations and relaxations are commonly used to improve computational tractability. However, we need to ensure that the approximations remain physically meaningful within the original problem. When addressing these challenges, there are trade-offs between the computational and solution performance. The focus of this paper is the second challenge.

In this work, we consider the interdependence of the power distribution network (PDN) and drinking water distribution network (WDN). Water supply pumps in the WDN are loads in the PDN. Around 4\% of the electricity consumption in the United States goes to pumping in drinking and waste water networks \cite{Denig-Chakroff2008}. However, this value can significantly vary over geographical location. For example, 19\% of energy consumption in California goes to water-related uses \cite{Klein2005}. Because most WDNs have water pumps at multiple locations and storage tanks that can store or supply water, we are able to shift pumping load spatially and temporally, and treat pumps as flexible loads. Here, we optimize and control water supply pumps to provide voltage support services to the PDN.

In the literature, ways of formulating and solving the integrated power-water optimization problem have been investigated in~\cite{ayyagari2021co, Zamzam2017, Fooladivanda2018,  OikonomouParvania2018, Liu2020}. 
In \cite{ayyagari2021co} and \cite{Zamzam2017}, the WDN and PDN are co-optimized to minimize power loss and distributed energy resource curtailment, respectively. Ref.~\cite{Fooladivanda2018} develops a framework for WDNs to respond to a signal to consume surplus renewable power. In~\cite{OikonomouParvania2018}, the authors use a group of WDNs to provide flexibility to the bulk transmission system. While power network modelling is not considered in \cite{Liu2020}, the authors calculate the demand response capacity of a WDN. Most work on integrated PDN-WDN optimization problems does not consider uncertainty. Uncertainty-aware optimization methods, such as robust and stochastic approaches, ensure that the systems are operating safely in the presence of uncertainty. In our prior work, we developed a  chance-constrained formulation that was solved via the scenario approach \cite{StuhlmacherProcIEEE} and an adjustable robust optimization formulation \cite{StuhlmacherCDC}. The scenario approach required a large amount of data and did not scale well to larger problems. The robust formulation was significantly more computationally tractable. However, a drawback to both of these approaches is that they tend to be excessively conservative. 

In this paper, we propose a computationally tractable and less conservative probabilistic approach that optimizes water pumping to provide voltage support under power demand uncertainty. Our approach combines a chance-constrained formulation to manage the uncertainty in the power distribution network and a probabilistically robust formulation to manage the impact of the PDN’s uncertainty on the WDN. Specifically, we develop an approach to schedule and control supply pump power to ensure that bus voltages in the PDN remain within their safe operating limits. The formulation ensures that the real-time control actions satisfy the network constraints at a user-specified probability level. To achieve computational tractability, we assume uncertainty distributions are known and build upon a robust approach for dissipative flow networks to analytically reformulate the probabilistic power and water constraints, resulting in a fast-to-solve deterministic problem. This enables application to large networks. The approach is much less conservative than robust approaches, but may perform poorly if assumed uncertainty distributions are inaccurate. We compare the performance (e.g., the computational time, cost, and empirical violation probabilities) of our approach with another uncertainty-aware approach and evaluate how practical these approaches would be when scaling the integrated PDN-WDN formulation to larger networks.

The contributions of this paper are the i) development of a tractable probabilistic formulation of the voltage support problem under power demand uncertainty subject to probabilistically robust WDN constraints and chance-constrained PDN constraints, 
ii) exploration of the performance and computational trade-offs between our proposed approach and other uncertainty-aware optimization methods, and iii) evaluation of the capability and practicality of the drinking water distribution network as a flexible load in a case study.

The remainder of the paper is organized as follows. Section~\ref{section: problem outline} provides an overview of the optimization framework and network modelling. In Section~\ref{section: formulations}, we present existing uncertainty-aware methods and develop our new probabilistic approach. We solve the uncertainty-aware approaches for a case study in Section~\ref{section: case study} and compare their performances. Lastly, we provide concluding remarks in Section~\ref{section: conclusion}.

\section{Problem Formulation} \label{section: problem outline}
Our goal is to optimize the water supply pump power consumption in a coupled PDN and WDN subject to network constraints and uncertain nodal power demands. We do not consider water demand uncertainty in this work. In~\cite{StuhlmacherProcIEEE}, we found that it is reasonable to assume that a portion of the tank capacity is pre-allocated to manage water demand uncertainty and therefore does not need to be explicitly considered in close-to-real-time operational planning problems. 

Our approach determines the scheduled water supply pump power consumption and the real-time pump power consumption adjustments needed to maintain safe voltage levels in the presence of power demand forecast errors, subject to the quasi-steady state power and water network constraints (i.e., steady-state operation within each time period of duration $\Delta T$). Specifically, we formulate an optimization problem that solves for the scheduled pump power consumption and affine control policy parameters that are used to adjust the pump power consumption in real-time as a function of the realization of the power demand forecast error vector. Using an affine control policy restricts the feasible space of real-time pump power consumption. However, an affine control policy eliminates the need for the water utility to re-solve the water pumping problem for each uncertainty realization and so it would be easier for the water utilities to implement. The optimization problem is of the form
\begin{mini!}|l|[2] 
{\mathbf{x}   }
%-
{  F(\mathbf{x}, \boldsymbol{\Delta\rho})} {}{} 
%-
\addConstraint{\mathcal{W}_1(\mathbf{x},\boldsymbol{\Delta \rho})}{\label{eqn: uncertain power constraints}}
%-
\addConstraint{\mathcal{W}_2(\mathbf{x},\boldsymbol{\Delta \rho}),}{\label{eqn: uncertain water constraints}}
\end{mini!}
%-
where $ F(\mathbf{x}, \boldsymbol{\Delta\rho})$ is the objective function. The constraint sets $\mathcal{W}_1(\mathbf{x}, \boldsymbol{\Delta\rho})$ and $\mathcal{W}_2(\mathbf{x},\boldsymbol{\Delta\rho})$ contain the quasi-steady state power flow and water flow constraints, respectively. These sets of constraints are functions of the decision variables $\mathbf{x}\in\mathbb{R}^d$ and the random variables $\boldsymbol{\Delta\rho}\in\mathbb{R}^n$. The following subsections detail these equations. 

Vector $\boldsymbol{\Delta \rho}$ contains the power demand forecast error $\Delta\rho_{k,\phi}^t$ for all buses~$k\in\mathcal{K}$ and phases~$\phi\in\Phi$ in the PDN and at all time periods within the scheduling horizon~$t \in \mathcal{T}$.
The power demand $\rho_{k,\phi}^t \in \boldsymbol{\rho}$ at bus~$k$, phase~$\phi$, and time~$t$ is the sum of the power demand forecast $\hat{\rho}_{k,\phi}^t \in \boldsymbol{\widehat{\rho}}$ and the power demand forecast error $\Delta\rho_{k,\phi}^t \in \boldsymbol{\Delta\rho}$.
Approaches to manage uncertainty in the optimization problem are presented in Section~\ref{section: formulations}. The benefits and trade-offs of each approach are then explored for a case study in Section~\ref{section: case study}.

\subsection{Power Distribution Network Constraints, $\mathcal{W}_1(\mathbf{x},\boldsymbol{\Delta \rho})$}\label{subsection: PDN constraints}
We consider radial, unbalanced, three-phase power distribution networks. In the PDN, we want to ensure that the voltages at each bus and phase are within their safe operating limits, i.e.,
\begin{align}
    (V_\text{min})^2\leq {Y}_{k,\phi}^t \leq (V_\text{max})^2\quad &\forall \,k\in \mathcal{K}, \phi\in\Phi, t\in\mathcal{T}, \label{eqn: voltage lim}
\end{align}
where ${Y}_{k,\phi}^t$ is the voltage magnitude squared at bus $k$, phase~$\phi$, and time~$t$. Parameters $V_\text{min}$ and $V_\text{max}$ are the lower and upper voltage magnitude limits. We use a linear power flow model to aid in the deterministic reformulations of the uncertainty-aware methods. In particular, a linear power flow model allows us to easily analytically reformulate the chance-constraints in Section~\ref{subsection: probabilistic approach}  and simplifies the robust reformulation of the power constraints in Section~\ref{subsection: robust}. The result is an equivalent set of deterministic constraints, assuming information on the uncertainty is known. We employ Lin3DistFlow \cite{Arnold, GanAndLow} -- a linearized, three-phase unbalanced power flow model -- to calculate the voltage magnitude squared
\begin{align}
    \boldsymbol{Y}_{k}^t&=\boldsymbol{Y}_{n}^t-\boldsymbol{M}_{kn}\boldsymbol{P}_{n}^t-\boldsymbol{N}_{kn}\boldsymbol{Q}_{n}^t &\forall \; k\in \mathcal{K}, t \in \mathcal T, \label{eqn: voltage}\\
    \boldsymbol{P}_{k}^t&=\boldsymbol{\rho}_k^t+\sum_{e\in\mathcal{P}_k}\boldsymbol{p}_{e}^t +\sum_{n\in\mathcal{I}_k}\boldsymbol{P}_{n}^t &\forall \; k\in \mathcal{K}, t \in \mathcal T, \label{eqn: real power flow}\\
    %-
    \boldsymbol{Q}_{k}^t&=\boldsymbol{\zeta}_k^t+ \sum_{e\in\mathcal{P}_k}\eta_{e}\boldsymbol{p}_e^t+\sum_{n\in\mathcal{I}_k}\boldsymbol{Q}_{n}^t &\forall \; k\in \mathcal{K}, t \in \mathcal T,  \label{eqn: reactive power flow}
\end{align} 
where $\boldsymbol{P}_{k}^t\in \mathbb{R}^3$ and $\boldsymbol{Q}_{k}^t\in\mathbb{R}^3$ are the three-phase real and reactive power flows entering bus~$k$ at time~$t$ and $\boldsymbol{Y}_k^t \in \mathbb{R}^3$ contains the three-phase voltage magnitude at bus~$k$ and time~$t$. Parameters $\boldsymbol{M}_{kn}\in\mathbb{R}^{3\times3}$ and $\boldsymbol{N}_{kn}\in\mathbb{R}^{3\times3}$ are calculated from the line impedance matrices for line $kn$. The set $\mathcal{I}_k\subseteq\mathcal{K}$ includes all buses that are directly downstream of bus~$k$. The three-phase real and reactive power demand at bus~$k$ and time~$t$ is denoted by $\boldsymbol{\rho}_k^t\in\mathbb{R}^{3}$ and $\boldsymbol{\zeta}_k^t\in\mathbb{R}^{3}$, respectively. The three-phase pump power consumption of pump~$e$ at time~$t$ is denoted by $\boldsymbol{p}_{e}^t\in\mathbb{R}^{3}$. The set $\mathcal{P}_k\subseteq \mathcal{P}$ contains all pumps that are connected to bus~$k$, where $\mathcal{P}$ contains all pumps in the WDN. We assume that the pumps are balanced loads and have a real-to-reactive power ratio $\eta_e$. It should be noted that other linear three-phase unbalanced power flow models could be used, e.g.,~\cite{Bernstein, Robbins2016}. Empirically, Lin3DistFlow and its variants perform very well. For example, \cite{VaninPSCC2020} found that the voltage magnitude accuracy of Lin3DistFlow consistently outperformed a second order cone relaxation in a case study of 500 existing Belgium feeders.

We utilize an affine control policy to adjust pump $e$'s single-phase pump power consumption $p_e^t$ in real time based on the power demand forecast error at time~$t$
\begin{align} 
    p_{e}^t = {p}_{\text{nom},e}^t + \boldsymbol{C}_{e}^t \boldsymbol{\Delta\rho}^t \quad \forall \, e \in \mathcal{P}, t\in \mathcal{T}, \label{eqn: voltage support control policy}
\end{align}
where decision variable $\boldsymbol{C}_e^t \in\mathbb{R}^{1\times n_t}$ is the control policy parameter row vector for pump $e$ at time~$t$ and $n_t$ is the number of random variables at time $t$. Decision variable ${p}_{\text{nom},e}^t$ is the single-phase scheduled pump power for pump $e$ at time~$t$. Note that the pump power's phase is not specified since we assume the pumps are balanced loads and so $p_{e}^t$ is the same in each phase. The control policy maps the change in power demand to a change in pump power consumption. 

The set of power constraints \eqref{eqn: uncertain power constraints} is composed of \eqref{eqn: voltage lim}-\eqref{eqn: voltage support control policy}. With the use of the affine control policy and Lin3DistFlow, the power constraints are linear. 

\subsection{Water Distribution Network Constraints, $\mathcal{W}_2(\mathbf{x},\boldsymbol{\Delta \rho})$}\label{subsection: WDN constraints}
The WDN can be represented as a connected directed graph composed of a set of pipes $\mathcal{E}$ and nodes $\mathcal{N}$. The pipes connect the nodes in the network, e.g., $ij\in\mathcal{E}$ is a pipe connecting node~$i$ to node~$j$. The water flow through a pipe may be positive or negative, where the sign indicates the direction that the water is moving. A pipe may contain a supply pump (i.e., $\mathcal{P}\subseteq\mathcal{E}$). A pump is restricted to a non-negative water flow rate since pumps can only pump in one direction. The nodes are composed of disjoint sets of reservoirs~$\mathcal{R}$, storage tanks~$\mathcal{S}$, and junctions~$\mathcal{J}$. We can characterize the operation of the WDN by the volumetric water flow rate through each pipe and the hydraulic head (which is equal to the sum of elevation and pressure head) at each node.

We approximate the head loss equations and the pump power curve to make the water constraints convex. While the robust approach discussed in Section~\ref{subsection: robust} and the proposed probabilistic approach in Section~\ref{subsection: probabilistic approach} do not require convex WDN constraints, we utilize the approximations to improve computation time. Additionally, the scenario-based probabilistic approach referenced in Section~\ref{subsection: probabilistic approach} requires convex constraints in order for the scenario approach to be applicable.  Refs.~\cite{StuhlmacherProcIEEE} and~\cite{Li2018} investigated the impact of the WDN approximations used and found that they were reasonable. For example,~\cite{StuhlmacherProcIEEE} empirically observed that the solutions of an approximated formulation using the linearized pump head-flow function \eqref{eqn: pump head gain} satisfied the original, nonconvex water constraints. 

To ensure safe operation, the hydraulic head at each node, the tank water levels, and pump flow rates must be within specified physical and/or operational limits at each time period~$t$
\begin{align}
    {H}_{\text{min},j}\leq {H}_j^t \leq {{H}}_{\text{max},j} &\quad \forall\, j\in\mathcal{N}, t\in\mathcal{T},  \label{eqn: head lim}\\
    {\ell}_{\text{min},j} \leq \ell_{j}^t \leq {\ell}_{\text{max},j} &\quad\forall\;j \in\mathcal{S}, t\in \mathcal T, \label{eqn: tank level bounds}\\
    0 \leq  {x}_{\text{min},ij} \leq x_{ij}^t \leq  {x}_{\text{max},ij} &\quad \forall\, ij\in\mathcal{P}, t\in \mathcal{T}, \label{eqn: flow limits}
\end{align} 
where $H_j^t$ is the hydraulic head at node $j$, $\ell_j^t$ is the water level in tank $j$, and $x_{ij}^t$ is the volumetric flow rate through pump~$ij$. The subscripts `min' and `max' on $H$, $\ell$, and $x$ indicate the lower and upper parameter limits of the variables, respectively. Additionally, we want to ensure that the tank levels are not depleted over the scheduling horizon, and so we require that the final tank levels are greater than or equal to the initial tank levels,  i.e.,
\begin{align}
    \ell_j^{t=|\mathcal{T}|} \geq \ell_j^{t=0} \quad \forall\, j\in\mathcal{S}. \label{eqn: final tank level}
\end{align}
The water flow can be determined from the quasi-steady state water flow equations \cite{Bhave} at each time period $t$
\begin{align}
    \sum_{i:ij\in \mathcal{E}} x_{ij}^t = -d_j^t  \quad &\forall\, j \in \mathcal{N}, t\in \mathcal{T},\label{eqn: conservation of water}\\
    \ell_{j}^{t} = \ell_{j}^{t-1} + \frac{\Delta T}{\gamma_j} \sum_{i:ij\in\mathcal{E}} x_{ij}^t \quad&\forall\;j \in\mathcal{S},  t\in \mathcal T, \label{eqn: tank level}\\
    H_j^t = \widehat{h}_j \quad &\forall\,j\in\mathcal{R}, t\in\mathcal{T}, \label{eqn: reservoir head}\\
    H_i^t - H_j^t = k_{ij} x_{ij}^t |x_{ij}^t| \quad & \forall \,ij\in \mathcal{E}\setminus\mathcal{P}, t\in \mathcal{T},\label{eqn: pipe head loss}\\
    H_j^t - H_i^t = m^1_{ij} x_{ij}^t + m_{ij}^0 \quad & \forall \,ij\in \mathcal{P}, t\in \mathcal{T},\label{eqn: pump head gain}\\
    p_{e}^t= g_{ij}^1 x_{ij}^t + g_{ij}^0 \quad &\forall\; e=ij\in\mathcal{P},  t\in \mathcal T, \label{eqn: pump power}
\end{align}
where $d_j^t$ is the water injection (consumer demand is nonpositive) at node~$j$, $\gamma_j$ is the cross-sectional area of tank~$j$,  $\widehat{h}_j$ is the elevation at node $j$, $k_{ij}$ is the resistance coefficient of pipe~$ij$, $m_{ij}^0$ and $m_{ij}^1$ are the head loss parameters of pump~$ij$, and $g_{ij}^1$ and $g_{ij}^0$ are parameters that relate power and water flow for pump~$ij$. The conservation of water is ensured in \eqref{eqn: conservation of water}. In \eqref{eqn: tank level}, the tank level is a function of the previous period's tank level and the tank's water injection in the current period.
Reservoirs are modelled as infinite sources and the hydraulic head is fixed in \eqref{eqn: reservoir head}.  The Darcy-Weisbach pipe head loss formulation, in which head loss is a function of the pipe flow rate, is given in 
\eqref{eqn: pipe head loss}. The pump head gain and single-phase pump power consumption as functions of the pump flow rate are given in~\eqref{eqn: pump head gain} and \eqref{eqn: pump power}, respectively. Here, we assume that the pump head and pump power are affine functions of water flow. 

In our formulation, we approximate the pipe head loss equation \eqref{eqn: pipe head loss} with a quasi-convex hull \cite{Li2018} in order to make the constraint convex. Equation \eqref{eqn: pipe head loss} is replaced with 
\begin{equation} \label{eqn: head loss pipe CONVEX}
\begin{aligned} 
H_i^t-H_j^t&\leq (2\sqrt{2}-2) k_{ij} {x}_{\text{max},ij}x_{ij}^t +(3-2\sqrt{2}) k_{ij} {x}_{\text{max},ij}^2, \\ 
H_i^t-H_j^t&\geq (2\sqrt{2}-2)k_{ij}|{x}_{\text{min},ij}|x_{ij}^t -(3-2\sqrt{2})k_{ij}{x}_{\text{min},ij}^2 ,\\
H_i^t-H_j^t&\geq 2k_{ij}{x}_{\text{max},ij}x_{ij}^t -k_{ij}{x}_{\text{max},ij}^2, \\ 
H_i^t-H_j^t&\leq 2k_{ij}|{x}_{\text{min},ij}|x_{ij}^t +k_{ij}{x}_{\text{min},ij}^2 ,
\end{aligned}
\end{equation}
where ${x}_{\text{min},ij}$ and ${x}_{\text{max},ij}$ are the lower and upper limits on pipe $ij$'s flow rate.

The set of water constraints \eqref{eqn: uncertain water constraints} is composed of \eqref{eqn: head lim}-\eqref{eqn: reservoir head}, \eqref{eqn: pump head gain}-\eqref{eqn: head loss pipe CONVEX}. 

\subsection{Objective Function, $ F(\mathbf{x}, \boldsymbol{\Delta\rho})$}\label{subsection: objective function}
The objective function minimizes the cost of the scheduled pump power and the cost associated with adjusting the pump power in real time (e.g., the wear-and-tear on the pumps from more frequent and larger magnitude real-time pump adjustments). In order to quantify the real-time pump adjustment cost, we define the range of real-time pump power adjustments due to the affine control policy \eqref{eqn: voltage support control policy}, i.e., the \textit{voltage support capacity constraints}
\begin{align}
    -\underline{R}_{e}^t\leq 3{\boldsymbol{C}}_{e}^t \boldsymbol{\Delta\rho}^t \leq \overline{R}_{e}^t\quad \forall \,e\in\mathcal{P},t\in\mathcal{T}, \label{eqn: R vs define}\\
    \underline{R}_{e}^t, \overline{R}_{e}^t\geq0\quad \forall \,e\in\mathcal{P},t\in\mathcal{T}, \label{eqn: R vs nonnegative}
\end{align}
where $\underline{{R}}_e^t, \overline{{R}}_e^t \in \mathbb{R}^{+}$ are the largest decrease and increase in pump $e$'s three-phase power consumption from the schedule. The voltage support capacity constraints are included with the reformulation of the robust power constraints in Section~\ref{subsection: robust} and the probabilistically robust water constraints in Section~\ref{subsection: probabilistic approach}. It is unclear how to best represent the real-time pump power cost; \cite{StuhlmacherProcIEEE} explored several options. However, if pump wear-and-tear is affected by the magnitude of the pump adjustments, then it is  reasonable to incorporate the real-time voltage support capacity into the cost function.

Our objective function can then be written as
\begin{align}
     F(\mathbf{x}, \boldsymbol{\Delta\rho}):= \sum_{t\in\mathcal{T}} \sum_{e\in\mathcal{P}} 3\pi^t  p_{\text{nom},e}^t + \pi_{\text{vs}}^t \left(\underline{R}_{e}^t  +\overline{R}_{e}^t \right)\label{eqn: obj}
\end{align}
where $\pi^t$ is the cost of electricity and $\pi_{\text{vs}}^t$ is the cost associated with the voltage support capacity at time~$t$.

\section{Approaches to Manage Uncertainty} \label{section: formulations}
Next, we discuss two uncertainty-aware methods to ensure that our solution will be feasible for a range of power demand forecast errors. First, we describe an existing adjustable robust formulation in Section~\ref{subsection: robust}.  Second, we develop our new probabilistic approach in Section~\ref{subsection: probabilistic approach}. Our approach builds on the robust reformulation from Section~\ref{subsection: robust} and existing approaches for chance-constraints. Third, Section~\ref{subsection: deterministic} describes a decoupled, deterministic approach, which provides a point of comparison for the uncertainty-aware methods. The benefits and disadvantages of these methods are explored in the case study in Section~\ref{section: case study}.

\subsection{Robust Approach} \label{subsection: robust}
We compare our probabilistic formulation with the adjustable robust approach of \cite{StuhlmacherCDC}. Furthermore, the monotonicity properties that are applied to the robust approach are also used to support the analytical reformulation of the probabilistically robust WDN constraints in the probabilistic approach in Section~\ref{subsection: probabilistic approach}.
The approach ensures that the power and water constraints are never violated under uncertainty, i.e., $\mathcal{W}_1(\mathbf{x},\boldsymbol{\Delta \rho})$ and $\mathcal{W}_2(\mathbf{x},\boldsymbol{\Delta \rho})$ are satisfied for all $\boldsymbol{\Delta \rho} \in \mathcal{U}$ where $\mathcal{U}$ is the uncertainty set. We assume the power demand forecast error is uncertain but bounded. 

For the power constraints, the power flow equalities \eqref{eqn: voltage}-\eqref{eqn: voltage support control policy} are substituted into the voltage limit inequalities \eqref{eqn: voltage lim} to form 
\begin{align}
    {V}_\text{min}^2 \leq {Y}_{k,\phi}^t(\mathbf{x}, \boldsymbol{\Delta\rho}) \leq V_\text{max}^2 \quad\forall\,k\in\mathcal{K},\phi\in\Phi,t\in\mathcal{T}.\label{eqn: y function}
\end{align}
The function ${Y}_{k,\phi}^t(\boldsymbol{x}, \boldsymbol{\Delta \rho})$ returns the voltage magnitude squared at bus~$k$, phase~$\phi$, and time~$t$ which is an affine function of the power demand forecast errors $\boldsymbol{\Delta \rho}$, the scheduled pump power $\boldsymbol{p}_\text{nom}$, and the control policy parameters $\boldsymbol{C}_e^t\in\boldsymbol{C}$ at time~$t$. 
The voltage support capacity constraints \eqref{eqn: R vs define}-\eqref{eqn: R vs nonnegative} and voltage limit constraints \eqref{eqn: y function} are linear in the random variables so we use explicit maximization \cite{lofberg2012} to obtain the robust counterpart. 

For the water constraints, we utilize monotonicity properties~\cite{Vuffray2015} to tractably reformulate the robust water constraints. This property holds in the WDN under several assumptions. First, we assume that the tank head is not strictly dependent on the tank level (e.g., there is a booster pump and/or valve connected to the tank). Second, we operate the WDN such that a positive deviation in all reservoir injections will cause a negative change in all tank injections. Third, the head loss and pump power consumption are increasing in flow rate. This last assumption requires that a pump's on/off status does not change in real time. The implications of these assumptions are discussed further in \cite{StuhlmacherCDC}. The monotonicity properties presented in~\cite{StuhlmacherCDC} state that if there is an increase in supply pump power consumption for all pumps, then the hydraulic head at all junctions can only decrease and the water levels at all tanks can only increase. Conversely, the hydraulic head can only increase and the tank levels can only decrease if there is a reduction in supply pump power consumption. From the monotonicity property, we only need to consider the feasibility of the extreme pump power cases to guarantee feasibility for all intermediary pump power set points. In \eqref{eqn: R vs define}, we defined the largest pump power adjustments, i.e., $\underline{R}_e^t \in \boldsymbol{\underline{R}}$ and $\overline{R}_e^t \in \boldsymbol{\overline{R}} \,\, \forall \,e\in\mathcal{P},\, t\in\mathcal{T}$. Therefore we can tractably rewrite the robust water constraints \eqref{eqn: uncertain water constraints} which hold for all uncertainty realizations $\boldsymbol{\Delta\rho}\in\mathcal{U}$ as 
\begin{subequations} \label{eqn: WF extreme sets}
\begin{align}
    \Gamma_\text{nom}({\boldsymbol{p}}_\text{nom})\label{eqn: scheduled water},\\
    \Gamma_\text{extreme}({\boldsymbol{p}}_\text{nom}+\boldsymbol{\overline{R}}), \label{eqn: extreme water max}\\
    \Gamma_\text{extreme}({\boldsymbol{p}}_\text{nom}-\boldsymbol{\underline{R}}), \label{eqn: extreme water min}%-
\end{align}
\end{subequations}
where $\Gamma_\text{nom}(\cdot)$  and $\Gamma_\text{extreme}(\cdot)$ are the set of WDN equations for the scheduled operation (i.e., \eqref{eqn: head lim}-\eqref{eqn: reservoir head}, \eqref{eqn: pump head gain}-\eqref{eqn: head loss pipe CONVEX}) and extreme operating conditions (i.e., \eqref{eqn: head lim}-\eqref{eqn: flow limits}, \eqref{eqn: conservation of water}-\eqref{eqn: reservoir head}, \eqref{eqn: pump head gain}-\eqref{eqn: head loss pipe CONVEX}). We use an overbar and underbar notation for the sets of decision variables associated with \eqref{eqn: extreme water max} and \eqref{eqn: extreme water min}, e.g., $\overline{\boldsymbol{H}}$ and $\underline{\boldsymbol{H}}$. 
Additional details are provided in \cite{StuhlmacherCDC}.

Then, the robust formulation is
\begin{mini*}|l|[2] 
{\mathbf{x}   }
%-
{ \eqref{eqn: obj} } {}{\tag*{\textbf{(R)}}} 
%-
%-
\addConstraint{\nu(\mathbf{x})\leq0}{}
\addConstraint{\eqref{eqn: WF extreme sets},}{}
\end{mini*}
where $\nu(\mathbf{x})\leq0$ represents the robust reformulation of the voltage support capacity and power constraints, i.e.,  \eqref{eqn: R vs define}-\eqref{eqn: R vs nonnegative}, \eqref{eqn: y function}. The decision variable $\mathbf{x}$ contains the scheduled pump power $\boldsymbol{p}_\text{nom}$; voltage support capacity $\boldsymbol{\underline{R}}$ and  
$\boldsymbol{\overline{R}}$; control policy parameters $\boldsymbol{C}$; hydraulic heads $H_j^t \in \boldsymbol{H}$, $\underline{H}_j^t \in \boldsymbol{\underline{H}}$, and $\overline{H}_j^t \in \boldsymbol{\overline{H}}$; tank levels 
$\ell_j^t \in \boldsymbol{\ell}$,  $\underline{\ell}_j^t \in \boldsymbol{\underline{\ell}}$, and $\overline{\ell}_j^t \in \boldsymbol{\overline{\ell}}$; and flow rates $x_{ij}^t \in \boldsymbol{x}$,  $\underline{x}_{ij}^t \in \boldsymbol{\underline{x}}$, and $\overline{x}_{ij}^t \in \boldsymbol{\overline{x}}$.

\subsection{Probabilistic Approach} \label{subsection: probabilistic approach}
We also consider a probabilistic approach. Here, we want to ensure that the network constraints are satisfied at high probability levels.

In previous work, the scenario approach \cite{Campi2009} was employed to solve a probabilistic formulation with convex network constraints \cite{StuhlmacherProcIEEE}. While there are certain advantages to the scenario approach (e.g., it does not require knowledge of the uncertainty distribution and enforces the chance-constraint jointly), it requires a large number of scenarios, which depends on the number of decision variables, and is not scalable to large networks and scheduling horizons. For example, we found that the scenario-based probabilistic approach was unable to solve the water pumping problem on a PDN-WDN comparable to that of the case study (see Fig.~\ref{Fig: wdn-pdn}) over a twelve-hour scheduling horizon due to memory issues~\cite{StuhlmacherPSCC}.

To improve computational tractability, we develop an uncertainty-aware optimization framework that contains chance-constrained power constraints and probabilistically robust water constraints.  Specifically, the chance-constrained power constraints guarantee feasibility of the PDN operation under power demand uncertainty with a high probability. The probabilistically robust water constraints guarantee feasibility of the real-time voltage support control policy on WDN operation with a high probability. Using probabilistic constraints is a reasonable approach since there are other components in the PDN that can provide voltage support. Additionally, small deviations outside of system limits for short periods of time may be acceptable.

We derive an analytical reformulation of the integrated PDN-WDN problem utilizing the monotonicity properties discussed in Section~\ref{subsection: robust} and assuming knowledge of the uncertainty distribution.
We rewrite the joint chance constraints as individual chance constraints. Unlike joint chance constraints, individual chance constraints do not significantly increase the computational complexity. Additionally, individual chance constraints are effective at reducing the joint violation probability and are less conservative than explicitly formulated joint chance constraints \cite{halilbavsic2018convex}. Furthermore, there is added flexibility in identifying active constraints and tuning constraint violation levels individually. We evaluate the empirical joint and individual reliability in the case study in Section~\ref{section: case study}.   

To analytically reformulate the power constraints, the voltage limit constraints are separated to form $2\times|\mathcal{K}|\times|\Phi|\times |\mathcal{T}|$ individual chance constraints 
\begin{align}
    \mathbb{P}\!\begin{bmatrix}
         {Y}_{k,\phi}^t(\mathbf{x}, \boldsymbol{\Delta\rho}) \leq {V}_\text{max}^2 \end{bmatrix}\! \geq 1-\epsilon_{\text{p}}\label{eqn:cc ymax} \,\,& \forall \, k\in\mathcal{K},\phi\in\Phi, t\in\mathcal{T},\\
    \mathbb{P}\!\begin{bmatrix}
        {V}_\text{min}^2 \leq {Y}_{k,\phi}^t(\mathbf{x}, \boldsymbol{\Delta\rho})  \end{bmatrix}\! \geq 1-\epsilon_{\text{p}} \label{eqn: cc ymin}\,\,&\forall \, k\in\mathcal{K},\phi\in\Phi, t\in\mathcal{T}.
\end{align}
The chance constraints ensure that each voltage limit is satisfied for a user-specified probability level $1-\epsilon_{\text{p}}$, where~$\epsilon_{\text{p}}$ is the individual violation level. It should be noted that the individual violation level can be different for each individual chance constraint. We analytically reformulate the power chance-constraints assuming that the power demand forecast error follows a normal distribution with a known mean $\mu \in\mathbb{R}^n$ and covariance $\Sigma \in \mathbb{R}^{n\times n}$. 
Ref.~\cite{roald2015security} found that a reformulation based on a normal distribution is reasonable (i.e., provides good trade-offs between cost and security) in systems with a large number of uncertain variables, e.g., uncertain power demands. The individual chance constraints can be equivalently written in the following deterministic form \cite{roald2015security, boyd} 
\begin{align}
    \mathbb{P} [  a(\mathbf{x}) + &b(\mathbf{x}) \boldsymbol{\Delta\rho} \leq c  ] \geq 1- \epsilon \label{eqn: cc individual}\\ 
    \Leftrightarrow \, &a(\mathbf{x}) \leq c-b(\mathbf{x}) \mu - f^{-1}(1-\epsilon) \| b(\mathbf{x})\Sigma^{1/2}\|_2,\notag
\end{align}
where $a(\mathbf{x})\in\mathbb{R}$ and $b(\mathbf{x})\in\mathbb{R}^{1\times n}$ are affine functions of the decision variables, $c\in\mathbb{R}$ is a constant, and $f^{-1}(\cdot)$ is either the inverse cumulative distribution function or a probability inequality \cite{roald2015security}. Since we assume a normal distribution, $f^{-1}(\cdot)$ is the inverse cumulative distribution function of the standard normal distribution.

For the water constraints \eqref{eqn: uncertain water constraints}, we build on a  probabilistically robust method presented in~\cite{Margellos2014} to make the probabilistically robust method applicable and tractable for our formulation. In~\cite{Margellos2014}, the authors determine an uncertainty set and solve for the probabilistically robust constraints in two sequential steps. Our formulation concurrently solves these two steps, which allows us to incorporate the probabilistically robust water constraints into an optimization framework that manages power constraints in a different way. In \cite{Margellos2014}, the uncertainty set is determined via a scenario-based approach, where the number of samples generated depends on the number of uncertainty sources and a user-specified violation level. Since our problem considers many sources of uncertainty (i.e., every bus and phase over the entire scheduling horizon), this approach would be very conservative for our formulation. It should also be noted that scenario-based approaches are generally conservative in practice \cite{Vrakopoulou2013}.  However, our approach may have more requirements on our formulation (for example, distribution assumptions) than the method in \cite{Margellos2014}.

To tractably reformulate the water constraints, we analytically solve for a bounded set $\mathcal{D}$ given the uncertainty $\boldsymbol{\Delta\rho}$. We ensure that  $\mathcal{D}$ encloses a user-specified probability density. The water constraints are then solved robustly, i.e., \eqref{eqn: WF extreme sets}, given set $\mathcal{D}$. Since the robust reformulation of the water constraints relies on the extreme pump powers, we solve for a range of extreme pump powers, i.e., $p_{\text{nom},e}^t - {\underline{R}}_e^t$ and $p_{\text{nom},e}^t + {\overline{R}}_e^t$ that need to be feasible in \eqref{eqn: extreme water max} and \eqref{eqn: extreme water min}. In other words, we ensure that the probability density of the real-time pump adjustments within $[\boldsymbol{\underline{R}}, \boldsymbol{\overline{R}}]$  
is at least $1-\epsilon_{\text{w}}$, i.e.,
\begin{align}
    \mathbb{P}& \begin{bmatrix}
        {\boldsymbol{C}_e^t} \boldsymbol{\Delta\rho}^t \leq \overline{{R}}_e^t\end{bmatrix} \geq 1-\epsilon_{\text{w}} &\quad \forall \, e\in\mathcal{P}, t\in\mathcal{T},\label{eqn: cc Rmax}\\
    \mathbb{P}& \begin{bmatrix}
        \text{-}\underline{{R}}_e^t\leq {\boldsymbol{C}_e^t} \boldsymbol{\Delta\rho}^t \end{bmatrix} \geq 1-\epsilon_{\text{w}} &\quad \forall \, e\in\mathcal{P}, t\in\mathcal{T},\label{eqn: cc Rmin}
\end{align}
where $\epsilon_{\text{w}}$ is the individual violation level. 
Therefore, the uncertainty set $\mathcal{D}$ is the hyper-rectangle of the voltage support capacities, i.e., $\mathcal{D} := \times_{e=1}^{|\mathcal{P}|} \times_{t=1}^{|\mathcal{T}|} [{\underline{R}}_e^t, {\overline{R}}_e^t]$ where $\times$ is the Cartesian product. 
Given the robust water constraints \eqref{eqn: WF extreme sets} are satisfied, the water flow constraints are feasible for pump power set points inside  $[\boldsymbol{p}_\text{nom}-\boldsymbol{\underline{R}},  \boldsymbol{p}_\text{nom} + \boldsymbol{\overline{R}}]$. 
Equations~\eqref{eqn: cc Rmax} and \eqref{eqn: cc Rmin} can be reformulated as deterministic constraints using the format presented in~\eqref{eqn: cc individual}.

The deterministic reformulation of the probabilistically robust water distribution network constraints and  the analytical reformulation of the voltage limit chance constraints are combined into a single optimization formulation
\begin{mini*}|l|[2] 
{\mathbf{x}   }
%-
{  \eqref{eqn: obj}} {}{\tag*{\textbf{(P)}}} 
%-
\addConstraint{\nu_\text{p}(\mathbf{x})\leq0}{}
\addConstraint{\nu_\text{w}(\mathbf{x})\leq0}{}
%-
\addConstraint{\eqref{eqn: WF extreme sets},}{}
\end{mini*}
%-
where $\nu_\text{p}(\mathbf{x})\leq0$ represents the analytical reformulation \eqref{eqn: cc individual} of the probabilistic voltage limit constraints \eqref{eqn:cc ymax}-\eqref{eqn: cc ymin} and $\nu_\text{w}(\mathbf{x})\leq0$ represents the analytical reformulation of the probabilistically robust  voltage support capacity constraints \eqref{eqn: cc Rmax}-\eqref{eqn: cc Rmin}.
The decision variable $\mathbf{x}$ contains $\boldsymbol{p}_\text{nom}$,  $\boldsymbol{\underline{R}}$,  $\boldsymbol{\overline{R}}$,  $\boldsymbol{C}$, $\boldsymbol{H}$,  $\boldsymbol{\underline{H}}$,  $\boldsymbol{\overline{H}}$,  $\boldsymbol{\ell}$,  $\boldsymbol{\underline{\ell}}$, $\boldsymbol{\overline{\ell}}$,  $\boldsymbol{x}$,     $\boldsymbol{\underline{x}}$,    and $\boldsymbol{\overline{x}}$. 

\subsection{Deterministic Approach}\label{subsection: deterministic}
Last, we present the decoupled formulation in which the WDN solves for the pump schedule with no knowledge of the PDN and the power demand uncertainty. We use the deterministic problem -- a formulation in which the PDN is not explicitly considered -- as a way to evaluate and compare the solutions and performance of the uncertainty-aware approaches presented in the previous subsections. Specifically, we evaluate the reduction in violation probabilities and the impact of the PDN on the water pumping. The deterministic, decoupled water pumping formulation is 
\begin{mini*}|l|[2] 
{\mathbf{x}   }
%-
{ \sum_{t\in\mathcal{T}} \sum_{e\in\mathcal{P}} 3\pi^t  p_{e}^t } {}{\tag*{\textbf{(D)}}} 
%-
%-
\addConstraint{\eqref{eqn: head lim}-\eqref{eqn: reservoir head},  \eqref{eqn: pump head gain}-\eqref{eqn: head loss pipe CONVEX},}{}
\end{mini*}
where $\mathbf{x}:=\{{p}_e^t \in \boldsymbol{p}, {H}_j^t \in \boldsymbol{H}, \ell_j^t \in \boldsymbol{\ell}, x_{ij}^t \in \boldsymbol{x}\}$. The pump power decision variable $p_e^t$ is equivalent to the scheduled pump power $p_{\text{nom},e}^t$ since there are no real-time pump adjustments.

\section{Case Study}\label{section: case study}
\begin{figure}[t]
\centering
\includegraphics[width=3in]{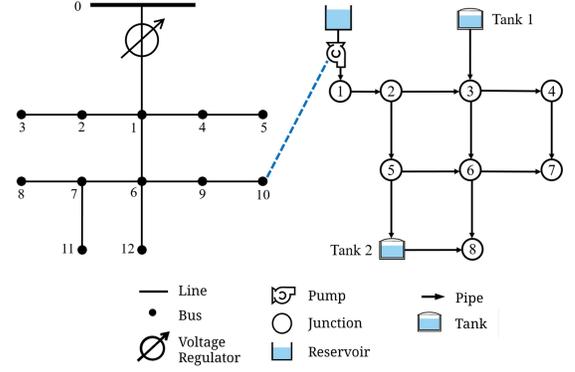}
\caption{Topology of the case study's coupled three-phase unbalanced PDN (left) and WDN (right). The blue dashed line indicates where the three-phase balanced water supply pump is located in the PDN.}
\label{Fig: wdn-pdn}
\end{figure}

In our case study, we consider the coupled PDN and WDN depicted in Fig.~\ref{Fig: wdn-pdn}.  We first describe the case study and then present the results. 

\subsection{Set Up} \label{subsection: setup}
\begin{figure}[t]
\centering
\includegraphics[width=3in]{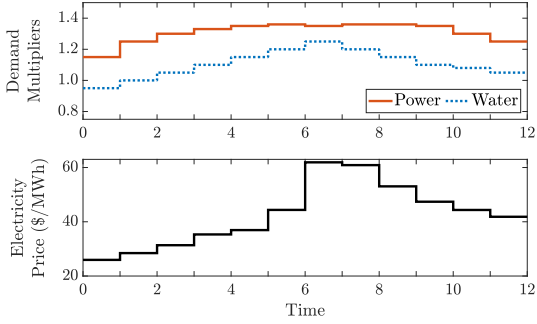}
\caption{Time-varying water and power demand multipliers (top) and electricity price (bottom) over the scheduling horizon.}
\label{Fig: Set up parameters}
\end{figure}

We use the coupled PDN and WDN from \cite{StuhlmacherCDC}. The WDN is based on an example network provided in EPANET, an open-source hydraulic modelling and simulation software program~\cite{EPANETGithub}. For the PDN, we use the IEEE 13-bus topology~\cite{IEEE13}. All modifications and parameter values are provided in~\cite{StuhlmacherCDC} except the following change to the water and power demands. In order to make the case study more realistic, the nominal water demand at each junction is multiplied by a time-varying constant. Similarly, the nominal power demand at each bus and phase is multiplied by a time-varying constant. The dashed blue curve and the solid red curve in the top plot of Fig.~\ref{Fig: Set up parameters} depict the water and power demand multipliers over a twelve-hour scheduling horizon where each time period has a duration $\Delta T$ of one hour. Additionally, we pulled electricity prices $
\pi^t$ for the Midcontinent Independent System Operator (MISO) for July 21st, 2021, 7:00-18:00 from \cite{MISO}. These electricity prices are shown in the bottom plot of Fig.~\ref{Fig: Set up parameters}. We set $\pi_\text{vs} = 5$~\$/MWh. The minimum and maximum voltage limits are 0.95~pu and 1.05~pu, respectively.

We consider two different types of uncertainty distributions -- a normal distribution and a student t-distribution. In the probabilistic approach (P) presented in Section~\ref{subsection: probabilistic approach}, we assume that the distribution of the power demand forecast error is normal, which may not be realistic. To observe how this approach performs when the actual distribution is not well known, we fit a multivariate normal distribution to 500 randomly drawn samples from the actual distribution, which is either a multivariate normal distribution or multivariate t-distribution. Student t-distributions have heavier tails than a normal distribution, leading to more extreme power demand forecast errors. It should be noted that we can also reformulate the chance-constraints assuming that the distribution is a t-distribution \cite{roald2015security}. However, a focus of this case study is evaluating the performance when the underlying uncertainty is different than what we assume. 

For the multivariate t-distribution, we set the degrees of freedom equal to 3, the correlation coefficients to 0.2, and scale the distribution such that the standard deviation is equal to $\alpha\hat{\rho}_{k,\phi}^t$ where $\alpha$ is a percentage and $\hat{\rho}_{k,\phi}^t$ is the forecasted power demand. Additionally, we truncate the t-distribution at ten times the standard deviation.
For the multivariate normal distribution, we generate samples from a load forecast error model that consists of a zero-mean normally distributed global error (with standard deviation $\sigma_0\hat{\rho}_{k,\phi}^t$) and a  zero-mean normally distributed nodal error (with standard deviation $\sigma_\text{node}\hat{\rho}_{k,\phi}^t$), similar to \cite{karangelos2020}. The normal distribution is truncated at three standard deviations from the mean. 
We consider three cases, which are described in Table~\ref{table: distributions}.

\begin{table}%[!t]
\renewcommand{\arraystretch}{1.3}
% if using array.sty, it might be a good idea to tweak the value of
% \extrarowheight as needed to properly center the text within the cells
\caption{Case Studies}
\label{table: distributions}
\centering
\begin{tabular}{ccc}
\hline
Case & Actual Distribution & Parameters\\
\hline
A & t-distribution & $\alpha=8$\%\\
B & t-distribution & $\alpha=1.5$\%\\
C & normal distribution & $\sigma_0=1.01$\%, $\sigma_\text{node}=3.98$\%\\
\hline
\end{tabular}
\end{table}

The estimated multivariate normal distribution used in the probabilistic approach (P) is fitted using maximum likelihood estimation from 500 randomly generated samples drawn from the actual distribution. We compare our solution to the adjustable robust approach (R), where the power demand is uncertain but bounded. The bounds are set to the maximum values of 2,000 samples randomly generated from the actual distribution, i.e., $\Delta \rho_{k,\phi}^t \in [- \Delta{\rho}_{\text{max},k,\phi}^t,  \Delta{\rho}_{\text{max},k,\phi}^t]$ for all $k\in \mathcal{K}$, $\phi\in\Phi$, and $t\in \mathcal{T}$ where $\Delta {\rho}^t_{\text{max},k,\phi} = \max_{i\in\mathcal{U}_{2000}} |\Delta \rho_{k,\phi,i}^t|$ and $\mathcal{U}_{2000}$ is the set of 2,000 randomly drawn samples.

To evaluate the solution performance of each approach, we use the Monte Carlo method to determine the empirical violation probabilities given the actual distribution and the fitted normal distribution. We calculate the empirical violation probability under both distributions so that we can evaluate the performance given the actual and expected distribution. We generate 50,000 power demand forecast error scenarios for the twelve-hour scheduling horizon and determine the corresponding real-time pump power adjustments from the pump schedule. We determine the joint empirical violation probability by determining whether there exists a water and/or power network violation for any scenario in the scenario set.

For simplicity, we set all individual violation levels for the voltage limits constraints~$\epsilon_\text{p}$ to be the same and all individual violation levels for the voltage support capacity constraints~$\epsilon_\text{w}$ to be the same. The problems are solved with the Gurobi solver~\cite{Gurobi} using the JuMP package in Julia on a computer with a 64-bit Intel i7 dual core CPU at 3.40 GHz and 16 GB RAM. 

\subsection{Results}\label{subsection: results}

\begin{table*}%[!t]
\renewcommand{\arraystretch}{1.3}
% if using array.sty, it might be a good idea to tweak the value of
% \extrarowheight as needed to properly center the text within the cells
\caption{Probabilistic and Robust Results}
\label{table: distribution 1 results}
\centering
\begin{tabular}{cccccrrrr}
\hline
Case& Problem & $\epsilon_\text{p}$ & $\epsilon_\text{w}$ & Solver Time& Total Cost& Scheduled Cost& Average $\boldsymbol{\overline{R}}$ & Average $\boldsymbol{\underline{R}}$ \\
& Type& (\%) & (\%) &(s) &  (\% Increase) & (\% Increase) & (kW) & (kW) \\
\hline
A& (P) & $1 \times 10^{\text{-}2}$ & $1 \times 10^{\text{-}2}$ & 1.36 & 11.53 & 4.90 & 184.90 & 183.37\\
& (P) & $5 \times 10^{\text{-}3}$ & $1 \times 10^{\text{-}3}$ & 1.46 & 15.40 & 7.63 & 216.77 & 215.13\\
& (P)   & $3 \times 10^{\text{-}3}$ & $1 \times 10^{\text{-}4}$ & 1.44 & 18.88 & 8.78 & 281.59 & 279.60\\ 
 & (R) & -- & -- & \multicolumn{5}{c}{---------------------------------Infeasible---------------------------------}\\
\hline
B& (P) & $1 \times 10^{\text{-}5}$ & $1 \times 10^{\text{-}5}$ & 0.49 & 0 & 0 & 0 & 0\\
 &(R) & -- & --& 0.40 & 2.86 & 1.54 & 36.80 & 36.80\\
\hline
C & (P) & $1 \times 10^{0}$ & $1 \times 10^{0}$& 0.48 & 0 & 0 & 0 & 0\\
 &(P) & $1 \times 10^{\text{-}3}$ & $1 \times 10^{\text{-}3}$& 0.94 & 0.49 & 0.49 & 0 & 0\\
&(P) & $1 \times 10^{\text{-}9}$ & $1 \times 10^{\text{-}9}$& 1.21 & 6.30 & 3.78 & 69.98 & 69.69\\
 &(R) & -- & --& 0.49 & 26.86 & 12.21 & 407.31 & 407.31\\
\hline
\end{tabular}
\end{table*}

We solve the robust (R), probabilistic (P), and deterministic~(D) formulations from Section~\ref{section: formulations}. Table~\ref{table: distribution 1 results} displays the solver time, the total cost \eqref{eqn: obj} and the cost of the pump schedule (i.e., the first term in \eqref{eqn: obj}) as percent increases from the deterministic cost (which has the same total and pump schedule cost), and the average three-phase up and down voltage support capacity over the twelve-hour scheduling horizon for the probabilistic and robust approaches as we vary the user-selected individual violation levels and the uncertainty distribution.
As the individual violation levels decrease and the uncertainty distribution's standard deviation increases, the pump operation shifts to a more expensive operating point and larger real-time control actions are needed to respond to the uncertainty realizations. As a result, the costs and the voltage support capacity increase. In Table~\ref{table: distribution 1 results}, Case~A has the largest forecast error standard deviation and costs while Case~B has the smallest forecast error standard deviation and costs. 
It should be noted that the real-time up and down voltage support capacities for the probabilistic approach are not equal because the fitted normal distribution is not necessarily zero mean.

\begin{figure}[t]
\centering
\includegraphics[width=3.3in]{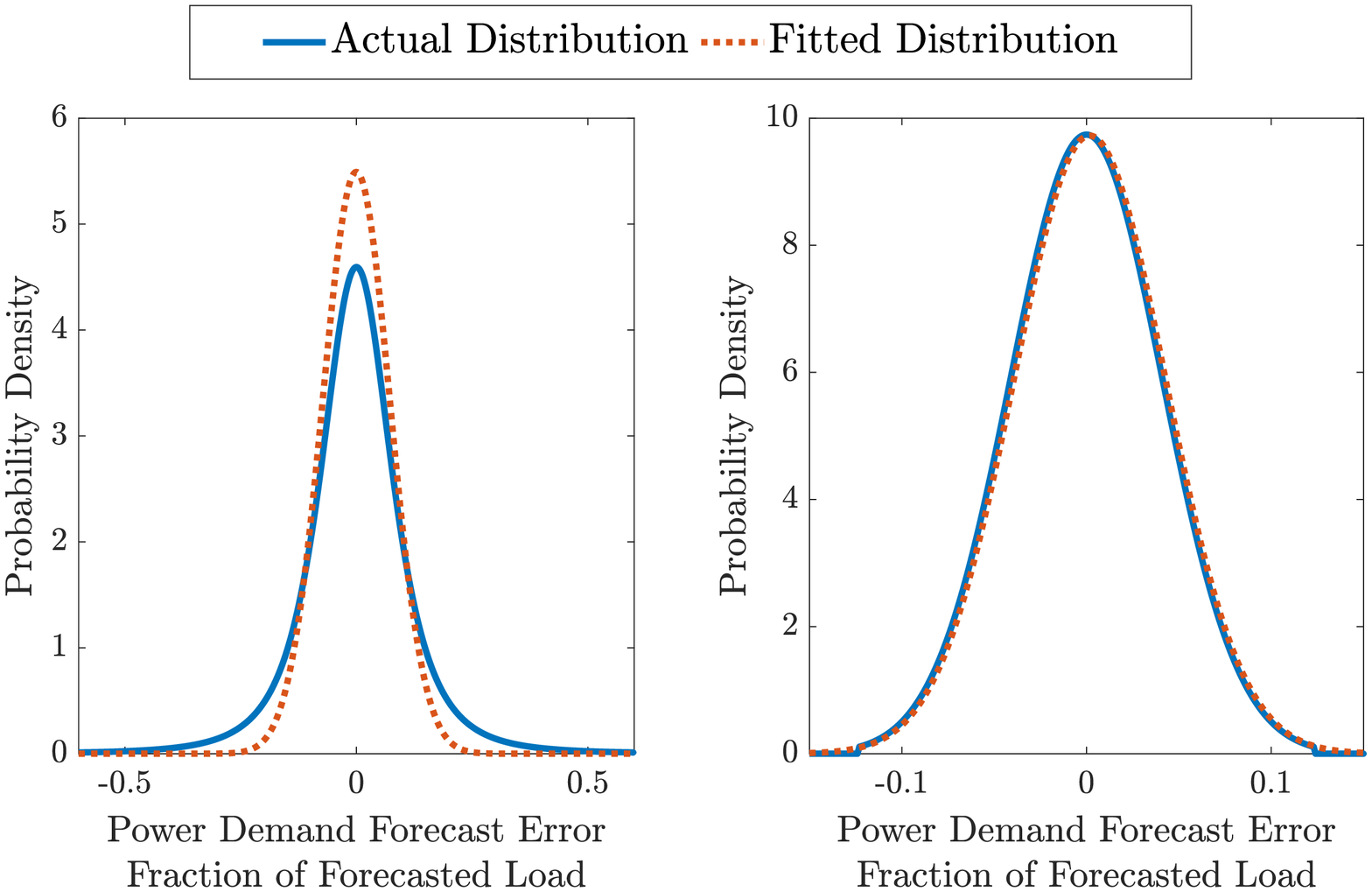}
\caption{Example probability density functions of the actual distributions (solid lines) for bus 8, phase $c$, and $t=1$ compared with the estimated normal distribution (dashed lines) fitted using maximum likelihood estimation with 500 samples. The actual distributions are Case~A (left) and Case~C (right).}
\label{Fig: pdfs}
\end{figure}

The robust solution requires a more expensive pump schedule and larger real-time voltage support capacity than the probabilistic solution. This is because the robust solution has to be feasible for all uncertainty realizations whereas the probabilistic solution needs to feasible for a certain probability density of the uncertainty realizations. For example, in Case~C, the total robust cost is a 26.86\% increase from the deterministic cost whereas the probabilistic case ($\epsilon_\text{p}$, $\epsilon_\text{w}=1 \times 10^{\text{-}9}$\%) is a 6.30\% increase from the deterministic cost. In many cases, the probabilistic formulation provides a solution whereas the robust formulation is infeasible. This illustrates the additional flexibility of a probabilistic solution versus a robust solution. However, the difference between the probabilistic and robust approaches intensifies when the fitted distribution is not representative of the actual distribution. This is illustrated in Fig.~\ref{Fig: pdfs}, which shows an example of the probability density functions for the actual distribution and its corresponding fitted normal distribution for Cases~A and~C. In Cases~A and B, the actual distributions have heavier tails than the fitted distributions. In these cases, the  probability of drawing a value from the fitted distribution that is near the robust bounds (which are formed from the actual distribution) is negligible. 
This indicates that the probabilistic approach using the fitted normal distribution may perform poorly if the assumed uncertainty distribution is inaccurate.

The solver times for the probabilistic and robust solutions are shown in Table~\ref{table: distribution 1 results}. The probabilistic and robust solver times are comparable. We find that the probabilistic approach would reasonably scale to larger networks since the probabilistic approach in the case studies solve in less than two seconds. For comparison, the scenario-based chance constrained formulation was unable to solve for a comparable PDN-WDN system due to memory issues.

\begin{figure}[t]
\centering
\includegraphics[width=\columnwidth]{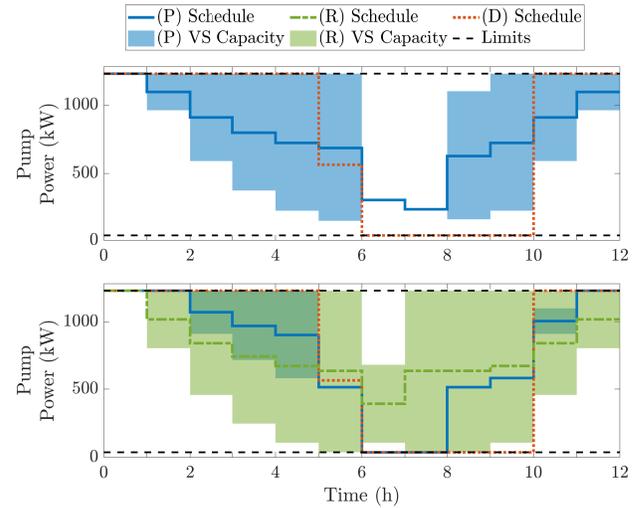}
\caption{Comparison of pump power schedules for the probabilistic (solid blue line), deterministic (dotted red lines), and robust (dashed green lines) solutions for (top) Case A ($\epsilon_\text{p}=3\times10^{\text{-}3}$\% and $\epsilon_\text{w}=1\times10^{\text{-}4}$\%) and (bottom) Case~C ( $\epsilon_\text{p}=1\times10^{\text{-}9}$\% and $\epsilon_\text{w}=1\times10^{\text{-}9}$\%). The blue and green bands indicate the range of real-time voltage support adjustments of the pump around the probabilistic and robust schedule, respectively. The black dashed lines indicate the pump power limits. }
\label{Fig: results plot}
\end{figure}

Fig.~\ref{Fig: results plot} illustrates how the uncertainty impacts the pump schedule for Case~A ($\epsilon_\text{p}=3\times10^{\text{-}3}$\% and $\epsilon_\text{w}=1\times10^{\text{-}4}$\%) and Case~C ($\epsilon_\text{p}$, $\epsilon_\text{w}=1 \times 10^{\text{-}9}$\%). We compare the probabilistic schedule (solid blue line), the deterministic schedule (dotted red line), and the robust schedule (dashed green line). The probabilistic and robust schedules vary less between time periods compared to the deterministic schedule and are more centered within the pump power limits. This is because the uncertainty-aware approaches need to respond to the uncertain power demand forecast error to ensure that the voltages remain within their limits. The range of real-time voltage support adjustments around the schedule are depicted with the blue and green bands for the probabilistic and robust solutions, respectively. In Case~A, the robust problem is infeasible and the pump in the probabilistic solution almost uses its full pumping range to provide real-time voltage support. In Case~C, the standard deviation of the forecast error is smaller so the probabilistic problem does not need as much voltage support capacity. Here, the robust solution is also feasible. As expected, the robust solution's real-time voltage support capacity is much larger than the  probabilistic solution's voltage support capacity.

\begin{table*}%[!t]
\renewcommand{\arraystretch}{1.3}
% if using array.sty, it might be a good idea to tweak the value of
% \extrarowheight as needed to properly center the text within the cells
\caption{Empirical Joint Violation Probabilities (\%)}
\label{table: empirical violation probability}
\centering
\begin{tabular}{ccccrrrrrr}
\hline
Case & Problem & $\epsilon_\text{p}$ & $\epsilon_\text{w}$ &\multicolumn{3}{c}{Fitted Distribution} & \multicolumn{3}{c}{Actual Distribution}\\
 & Type & (\%) & (\%)&Joint Power & Joint Water & Total & Joint Power & Joint Water&Total\\
\hline
A&(D) &  -- & -- &5.37 & 0&  5.37 & 11.11 & 0 & 11.11\\
&(P) &  $1 \times 10^{\text{-}2}$ & $1 \times 10^{\text{-}2}$ &0.04 & 0.01&  0.05 & 5.02 & 2.40 & 7.12\\
&(P) & $5 \times 10^{\text{-}3}$ & $1 \times 10^{\text{-}3}$ &0.02 & 0& 0.02 & 4.75 & 1.69 & 6.19\\
&(P) & $3 \times 10^{\text{-}3}$ & $1 \times 10^{\text{-}4}$ &0.01 & 0& 0.01 & 4.43 & 1.25 & 5.49\\
&(R) &  -- & -- &\multicolumn{6}{c}{----------------------------------Infeasible----------------------------------}\\
\hline
B&(D) & -- & -- &0 & 0& 0 & 0 & 0 & 0\\
&(P)& $1 \times 10^{\text{-}5}$ & $1 \times 10^{\text{-}5}$ &0 & 0& 0 & 0 & 0 & 0\\
&(R) & --& -- & 0 & 0 & 0 & 0 &0&0\\
\hline
C&(D) & -- & -- &0 & 0& 0 & 0.01 & 0 & 0.01\\
&(P) & $1 \times 10^{0}$ & $1 \times 10^{0}$ &0 & 0& 0 & 0.01 & 0 & 0.01\\
&(P) &  $1 \times 10^{\text{-}3}$ & $1 \times 10^{\text{-}3}$ &0 & 0& 0 & 0 & 0 & 0\\
&(P) &  $1 \times 10^{\text{-}9}$ & $1 \times 10^{\text{-}9}$ &0 & 0& 0 & 0 & 0 & 0\\
&(R) & -- & -- & 0 & 0 & 0 & 0 &0&0\\
\hline
\end{tabular}
\end{table*}

Next, we verify that the empirical violation probabilities of the individual chance constraints are below the user-selected violation levels, $\epsilon_\text{p}$ and $\epsilon_\text{w}$ when the samples are generated from the fitted, multivariate normal distribution. In general, we found that many empirical violation probabilities are well below the user-specified violation level. To demonstrate this, let's consider the voltage constraints in Case~A ($\epsilon_\text{p}=5\times10^{\text{-}3}$\% and $\epsilon_\text{w}=1\times10^{\text{-}3}$\%). The empirical voltage violation probabilities range from 0\% to $4\times10^{\text{-}3}$\%, with a mean of $8.8\times10^{\text{-}5}$\% and standard deviation of $4.9\times10^{\text{-}4}$\%. This indicates that many individual chance constraints are never active. We found that the violations only occurred on phase $c$ at the end of the network (buses 6-12), where a majority of all violations (around 47\%) occurred at bus~8. 

Table~\ref{table: empirical violation probability} evaluates the joint empirical violation probabilities. The joint water, joint power, and overall joint empirical violation probabilities are displayed for both the fitted and actual distributions. We observe that Cases B and C both have probabilistic solutions that have the same empirical violation probabilities as the robust solution, despite having lower costs. As expected, the empirical violation probability decreases as the violation level decreases. Additionally, the actual distribution generally had worse empirical violation probabilities than the fitted distribution when the fitted distribution was not representative of the actual distribution's probability distribution (e.g., Case~A). 
However, we observe that by decreasing the individual violation levels $\epsilon_\text{w}$ and $\epsilon_\text{p}$, we are able to improve the joint empirical violation probability of the actual distribution, despite the distribution being different than we expect. 
We also compare the empirical violation probabilities of the probabilistic approach with the  deterministic approach. The joint power and overall joint empirical violations are generally much higher in the deterministic approach since the WDN and PDN operation is completely decoupled. This illustrates the benefits of coupled operation to the power network and the coupled power-water system as a whole. As expected, the WDN does not experience any constraint violations in the deterministic case. In the probabilistic case, water constraint violations can occur since WDN operation is impacted by the uncertainty in the PDN. This can be seen in the joint water violation probabilities for Case~A in Table~\ref{table: empirical violation probability}.  
However, the empirical violation probabilities can be tuned by the individual constraint violation levels to provide a certain amount of reliability to each network. 
A benefit of individual chance constraints over joint chance constraints is that there is more flexibility in identifying and tuning important individual constraints for system security~\cite{halilbavsic2018convex}.

\section{Conclusion}\label{section: conclusion}
In this paper, we proposed a computationally tractable analytical reformulation of a probabilistic water pumping problem to provide voltage support to the PDN. The problem is subject to the probabilistically robust WDN constraints and chance-constrained PDN constraints and manages power demand uncertainty. We compared this approach with an adjustable robust method to evaluate the computational and solution performance. We found that the probabilistic approach is significantly less conservative than the robust approach and  that it is able to find solutions when the robust approach is infeasible. By adjusting individual violation levels, we can target network constraints that are important for system reliability. A drawback to the probabilistic approach is that we assume the uncertainty distribution is known, and so the probabilistic approach may perform poorly if the distribution is inaccurate. 
The probabilistic approach has a comparable computational performance to the robust approach. Therefore, this approach can be applied to large networks. In future work, we plan to apply the probabilistic formulation to larger networks and explore methods to minimize communication between the water and power system operators.  

%\section*{Acknowledgment}

\IEEEtriggeratref{17}


%\bibliographystyle{IEEEtran}

\end{document}